# Differential Transformation of Nonlinear Power Flow Equations

Yang Liu, *Student Member, IEEE*, Kai Sun, *Senior Member, IEEE*

***Abstract*— Differential transformation (DT) method has shown to be promising for power system simulation in our recent works. This letter applies the DT method to nonlinear power flow equations and proves that the nonlinear power flow equations are converted to formally linear equations after DT. This letter further extends the constant power load assumption in the classical power flow formulation to the more general and practical ZIP load assumption and proves that the linear equation after DT still holds for ZIP load model with slight modifications on the coefficient matrices. These results demonstrate great potential of the DT method in solving power system nonlinear problems.**

## I. Introduction

Differential transformation (DT) method is a promising approach for power system dynamic simulation. It is introduced to power system field in our recent work and has proved to be effective to solve both the ordinary differential equation (ODE) model [1] and the differential algebraic equation (DAE) model [2]. For a DAE model in transient stability simulation, the nonlinear current injection equations of synchronous generators, ZIP loads [2] and motor loads [3] were proved to satisfy a formally linear equation after DT, which is an important property and has enabled a highly efficient algorithm to solve the DAE model without iterations. This letter further applies the DTs to the nonlinear power flow equations with both constant power load model and ZIP load model. For both load models, we prove that the power flow equations for PQ buses, PV buses and the slack bus all satisfy formally linear equations after DT. We further derived the analytical expression of the coefficient matrices in the linear equation. These results demonstrate great potential of the DT method in solving power system nonlinear problems.

## II. Main Results

The Differential Transformation (DT) method was introduced to the power system field in [1]. A detailed introduction of the DT method and the transformation rules for various generic nonlinear functions in power system models are provided in [1]-[2].

### A. Nonlinear Power Flow Equations

The classical power flow equations [4]-[6] are given in (1a) where $\overline{\mathrm{S}}$ is a vector of the complex power injections, $\overline{\mathrm{V}}$ is a vector of bus voltage phasor, and $\mathbf{Y}_{\mathbf{bus}}$ is the bus admittance matrix. This letter considers a more general setting in (1b) which is obtained by adding the product of a loading parameter $\lambda$ and a constant vector $\overline{\mathrm{b}}$ to the left-hand side [6].

$$\begin{aligned} \overline{\mathrm{S}} &= \overline{\mathrm{V}}(\mathbf{Y}_{\mathrm{bus}}\overline{\mathrm{V}})^* && \text{(a)} \\ \overline{\mathrm{S}}+\lambda\overline{\mathrm{b}} &= \overline{\mathrm{V}}(\mathbf{Y}_{\mathrm{bus}}\overline{\mathrm{V}})^* && \text{(b)} \end{aligned} \qquad (1)$$

Equation (1b) is further written as the general form in (2) where $\boldsymbol{g}$ is a nonlinear vector field; $\boldsymbol{y}$ is the bus voltage vector under rectangular coordinates defined as $\boldsymbol{y}=[\boldsymbol{e}^{\mathrm{T}},\boldsymbol{f}^{\mathrm{T}}]^{\mathrm{T}}$, where $\boldsymbol{e}=[e_1,\dots, e_N]^{\mathrm{T}}$ and $\boldsymbol{f}=[f_1,\dots, f_N]^{\mathrm{T}}$ are respectively the real and imaginary parts of the bus voltage phasor; $N$ is the total number of buses; $\lambda$ is the loading parameter.

$$0 = \boldsymbol{g}(\boldsymbol{y},\lambda) \qquad (2)$$

The goal is to solve the power flow solution $\boldsymbol{y}$ under a given loading condition $\lambda$, or more generally, to determine how power flow solution $\boldsymbol{y}$ changes with loading parameter $\lambda$, shown in (3). Generally, analytical expression of (3) is unavailable due to the nonlinearity of $\boldsymbol{g}$ in (2), and iteration methods such as the family of Newton method are needed.

$$\boldsymbol{y} = \boldsymbol{y}(\lambda) \qquad (3)$$

### B. Overview of the Major Results

For the convenience, variables $\boldsymbol{y}$ and $\lambda$ in (2) can be interpreted as functions of time, i.e., $\boldsymbol{y}(t)$ and $\lambda(t)$ as shown in (4a), and the DTs of them are denoted by corresponding capital letters $\boldsymbol{Y}(k)$ and $\Lambda(k)$, where $k$ is the order of power series coefficients. This letter first derives DTs of (4a) with the general form in (4b) about the 0 to $k^{\mathrm{th}}$ order power series coefficients $\boldsymbol{Y}(0{:}k)$ and $\Lambda(0{:}k)$. Then, we further prove that (4b) satisfies a formally linear equation in (4c) where the matrices $\boldsymbol{A}$ and $\boldsymbol{B}$ only depends on 0 to $(k\text{-}1)^{\mathrm{th}}$ order power series coefficients $\boldsymbol{Y}(0{:}k\text{-}1)$ and $\Lambda(0{:}k\text{-}1)$. The derivation and proofs are conducted for constant power load assumption and ZIP load assumption respectively.

This work was supported in part by the ERC Program of the NSF and DOE under NSF Grant EEC-1041877 and in part by NSF Grant ECCS-1610025.

Y. Liu and K. Sun are with the Department of EECS, University of Tennessee, Knoxville, TN 37996 USA (e-mail: yliu161@vols.utk.edu, kaisun@utk.edu).

$$\begin{aligned} & 0 = \boldsymbol{g}(\boldsymbol{y}(t), \lambda(t)) && (a) \\ & \downarrow \\ & 0 = \boldsymbol{G}(\boldsymbol{Y}(0:k), \Lambda(0:k)) && (b) \\ & \downarrow \\ & 0 = \boldsymbol{A}_{gy}\boldsymbol{Y}(k) + \boldsymbol{A}_{g\lambda}\Lambda(k) + \boldsymbol{B}_g && (c) \end{aligned} \quad (4)$$

### *C. DTs of the Nonlinear Power Flow Equations with Constant Power Load Model*

Under the assumption of constant power load, the nonlinear power flow equation (2) is written into (5)-(8) in rectangular coordinates, where $\Omega_{PQ}$, $\Omega_{PV}$, $\Omega_{REF}$ are the set of PQ buses, PV buses and reference bus respectively, $p$ and $q$ are active and reactive power, $e$ and $f$ are the real and imaginary parts of bus voltages, $g$ and $b$ are real and imaginary parts of the admittance, $v$ is the voltage magnitude, superscript $sp$ means the value is specified, subscript $i$ and $j$ are the index of buses.

$$\begin{aligned} p_i^{sp} &= g_p(\boldsymbol{y}, \lambda) \\ &= -\lambda\Delta p_i + \sum_{j=1}^{N} g_{ij}\left(e_i e_j + f_i f_j\right) + \sum_{j=1}^{N} b_{ij}\left(f_i e_j - e_i f_j\right) \\ &\text{if } i \in \Omega_{PQ} \cup \Omega_{PV} \end{aligned} \quad (5)$$

$$\begin{aligned} q_i^{sp} &= g_q(\boldsymbol{y}, \lambda) \\ &= -\lambda\Delta q_i - \sum_{j=1}^{N} b_{ij}\left(e_i e_j + f_i f_j\right) + \sum_{j=1}^{N} g_{ij}\left(f_i e_j - e_i f_j\right) \\ &\text{if } i \in \Omega_{PQ} \end{aligned} \quad (6)$$

$$(v_i^{sp})^2 = g_v(\boldsymbol{y}) = e_i^2 + f_i^2, \text{if } i \in \Omega_{PV} \quad (7)$$

$$\begin{aligned} e_i^{sp} &= g_e(\boldsymbol{y}) = e_i, \text{if } i \in \Omega_{REF} \\ f_i^{sp} &= g_f(\boldsymbol{y}) = f_i, \text{if } i \in \Omega_{REF} \end{aligned} \quad (8)$$

Following the transformation rules in [1]-[2], the DTs of (5)-(8) are derived in (9)-(12) respectively with details omitted.

$$\begin{aligned} p_i^{sp}\rho(k) &= G_p(\boldsymbol{Y}, \Lambda) \\ &= -\Delta p_i \Lambda(k) + \sum_{j=1}^{N} g_{ij}\left(E_i(k) \otimes E_j(k) + F_i(k) \otimes F_j(k)\right) \\ &+ \sum_{j=1}^{N} b_{ij}\left(F_i(k) \otimes E_j(k) - E_i(k) \otimes F_j(k)\right) \\ &\text{if } i \in \Omega_{PQ} \cup \Omega_{PV} \end{aligned} \quad (9)$$

$$\begin{aligned} q_i^{sp}\rho(k) &= G_q(\boldsymbol{Y}, \Lambda) \\ &= -\Delta q_i \Lambda(k) - \sum_{j=1}^{N} b_{ij}\left(E_i(k) \otimes E_j(k) + F_i(k) \otimes F_j(k)\right) \\ &+ \sum_{j=1}^{N} g_{ij}\left(F_i(k) \otimes E_j(k) - E_i(k) \otimes F_j(k)\right) \\ &\text{if } i \in \Omega_{PQ} \end{aligned} \quad (10)$$

$$\begin{aligned} (v_i^{sp})^2\rho(k) &= G_v(\boldsymbol{Y}) \\ &= E_i(k) \otimes E_i(k) + F_i(k) \otimes F_i(k), \text{ if } i \in \Omega_{PV} \end{aligned} \quad (11)$$

$$\begin{aligned} e_i^{sp}\rho(k) &= G_e(\boldsymbol{Y}) = E_i(k), \text{if } i \in \Omega_{REF} \\ f_i^{sp}\rho(k) &= G_f(\boldsymbol{Y}) = F_i(k), \text{if } i \in \Omega_{REF} \end{aligned} \quad (12)$$

**Proposition 1**: The transformed power flow equations (9)-(12) respectively satisfy formally linear equations (13)-(16).

$$0 = \boldsymbol{a}_{\mathrm{P},i}\boldsymbol{Y}(k) - \Delta p_i \Lambda(k) + \varepsilon_i, \text{if } i \in \Omega_{PQ} \cup \Omega_{PV} \quad (13)$$

$$0 = \boldsymbol{a}_{\mathrm{Q},i}\boldsymbol{Y}(k) - \Delta q_i \Lambda(k) + \mu_i, \text{if } i \in \Omega_{PQ} \quad (14)$$

$$0 = \boldsymbol{a}_{\mathrm{V},i}\boldsymbol{Y}(k) + 0\Lambda(k) + \varsigma_i, \text{if } i \in \Omega_{PV} \quad (15)$$

$$\begin{aligned} 0 &= \boldsymbol{a}_{\mathrm{E},i}\boldsymbol{Y}(k) + 0\Lambda(k) - e_i^{sp}\rho(k), \text{if } i \in \Omega_{REF} \\ 0 &= \boldsymbol{a}_{\mathrm{F},i}\boldsymbol{Y}(k) + 0\Lambda(k) - f_i^{sp}\rho(k), \text{if } i \in \Omega_{REF} \end{aligned} \quad (16)$$

where $\boldsymbol{Y}(k) \in \mathbb{R}^{2N\times 1}$ and $\Lambda(k) \in \mathbb{R}$ are variables representing the DT of $\boldsymbol{y}$ and $\lambda$ respectively; $\boldsymbol{a}_{\mathrm{P},i}$, $\boldsymbol{a}_{\mathrm{Q},i}$, $\boldsymbol{a}_{\mathrm{V},i}$ $\boldsymbol{a}_{\mathrm{E},i}$, $\boldsymbol{a}_{\mathrm{F},i} \in \mathbb{R}^{1\times 2N}$ and $\varepsilon_i$, $\mu_i$, $\zeta_i \in \mathbb{R}$ are parameters given in the Appendix, together with the detailed proofs.

**Proposition 2**: Under the assumption of constant power load, the DTs (4b) of nonlinear power flow equations (4a) satisfy a formally linear equation (4c) with matrices $\boldsymbol{A}_{gy}$, $\mathbf{A}_{g\lambda}$, and $\boldsymbol{B}_g$ given by (17).

$$\boldsymbol{A}_{gy} = \begin{bmatrix} \boldsymbol{A}_{y,\mathrm{PQ}} \\ \boldsymbol{A}_{y,\mathrm{PV}} \\ \boldsymbol{A}_{y,\mathrm{REF}} \end{bmatrix}, \boldsymbol{A}_{g\lambda} = \begin{bmatrix} \boldsymbol{A}_{\lambda,\mathrm{PQ}} \\ \boldsymbol{A}_{\lambda,\mathrm{PV}} \\ \boldsymbol{A}_{\lambda,\mathrm{REF}} \end{bmatrix}, \boldsymbol{B}_g = \begin{bmatrix} \boldsymbol{B}_{\mathrm{PQ}} \\ \boldsymbol{B}_{\mathrm{PV}} \\ \boldsymbol{B}_{\mathrm{REF}} \end{bmatrix}$$

$$where$$

$$\boldsymbol{A}_{y,\mathrm{PQ}} = \begin{bmatrix} \boldsymbol{a}_{\mathrm{P},1} \\ \boldsymbol{a}_{\mathrm{Q},1} \\ \vdots \\ \boldsymbol{a}_{\mathrm{P},M} \\ \boldsymbol{a}_{\mathrm{Q},M} \end{bmatrix} \quad \boldsymbol{A}_{\lambda,\mathrm{PQ}} = -\begin{bmatrix} \Delta p_1 \\ \Delta q_1 \\ \vdots \\ \Delta p_M \\ \Delta q_M \end{bmatrix} \quad \boldsymbol{B}_{\mathrm{PQ}} = \begin{bmatrix} \varepsilon_1 \\ \mu_1 \\ \vdots \\ \varepsilon_M \\ \mu_M \end{bmatrix}$$

$$\boldsymbol{A}_{y,\mathrm{PV}} = \begin{bmatrix} \boldsymbol{a}_{\mathrm{P},M+1} \\ \boldsymbol{a}_{\mathrm{V},M+1} \\ \vdots \\ \boldsymbol{a}_{\mathrm{P},N-1} \\ \boldsymbol{a}_{\mathrm{V},N-1} \end{bmatrix} \quad \boldsymbol{A}_{\lambda,\mathrm{PV}} = -\begin{bmatrix} \Delta p_{M+1} \\ 0 \\ \vdots \\ \Delta p_{N-1} \\ 0 \end{bmatrix} \quad \boldsymbol{B}_{\mathrm{PV}} = \begin{bmatrix} \varepsilon_{M+1} \\ \zeta_{M+1} \\ \vdots \\ \varepsilon_{N-1} \\ \zeta_{N-1} \end{bmatrix}$$

$$\boldsymbol{A}_{y,\mathrm{REF}} = \begin{bmatrix} \boldsymbol{a}_{\mathrm{E},N} \\ \boldsymbol{a}_{\mathrm{F},N} \end{bmatrix} \quad \boldsymbol{A}_{\lambda,\mathrm{REF}} = \begin{bmatrix} \boldsymbol{0} \\ \boldsymbol{0} \end{bmatrix} \quad \boldsymbol{B}_{\mathrm{REF}} = \begin{bmatrix} -e_N^{sp}\rho(k) \\ -f_N^{sp}\rho(k) \end{bmatrix} \quad (17)$$

For notation simplicity, here we let buses 1 to $M$ be PQ buses, buses $M$+1 to $N$-1 be PV buses and bus $N$ be the reference bus. The Proposition 2 can be easily proved from Proposition 1 and the details are omitted.

### *D. DTs of the Nonlinear Power Flow Equations with ZIP Load Model*

Using the ZIP load model assumption, the power flow equations (5)-(8) needs to be modified. Specifically, (7)-(8) remain unchanged, but (5)-(6) need to be rewritten as (18)-(20) respectively, where $p_{i,ZIP}(\boldsymbol{y})$ and $q_{i,ZIP}(\boldsymbol{y})$ are derived in the subsection II-D-1) below. Note that the power injection in a ZIP load model is a function of the bus voltage vector ***y***, which is different from the constant power load where the power injection of loads are constants. To sum up, the power flow equations for ZIP load model are given by (18)-(20) and (7)-(8).

$$p_i^{sp} = g_p(\boldsymbol{y},\lambda) = -\lambda\Delta p_i + \sum_{j=1}^{N} g_{ij}\left(e_ie_j + f_if_j\right) + \sum_{j=1}^{N} b_{ij}\left(f_ie_j - e_if_j\right) \quad \text{if } i\in\Omega_{PV} \tag{18}$$

$$p_{i,ZIP}(\boldsymbol{y}) = g_p(\boldsymbol{y},\lambda) = -\lambda\Delta p_i + \sum_{j=1}^{N} g_{ij}\left(e_ie_j + f_if_j\right) + \sum_{j=1}^{N} b_{ij}\left(f_ie_j - e_if_j\right) \quad \text{if } i\in\Omega_{PQ} \tag{19}$$

$$q_{i,ZIP}(\boldsymbol{y}) = g_q(\boldsymbol{y},\lambda) = -\lambda\Delta q_i - \sum_{j=1}^{N} b_{ij}\left(e_ie_j + f_if_j\right) + \sum_{j=1}^{N} g_{ij}\left(f_ie_j - e_if_j\right) \quad \text{if } i\in\Omega_{PQ} \tag{20}$$

DTs of power flow equations with ZIP load model can be derived as follows. Compared with the power flow equation with constant power load model, the only difference lies in the left-hand-side in (19)-(20), i.e., $p_{i,ZIP}(\boldsymbol{y})$ and $q_{i,ZIP}(\boldsymbol{y})$. Therefore, we first derive DTs of $p_{i,ZIP}(\boldsymbol{y})$ and $q_{i,ZIP}(\boldsymbol{y})$ (denoted by $P_{i,ZIP}(\boldsymbol{Y})$ and $Q_{i,ZIP}(\boldsymbol{Y})$) in the following subsection II-D-2), then we will present the DTs of the full set of power flow equations (18)-(20) and (7)-(8) for ZIP load model in subsection II-D-3).

*1) Power Injection of ZIP Load Model*

For constant impedance component, the impedance $z_i^{sp} = \mathrm{Re}(z_i^{sp}) + j\,\mathrm{Im}(z_i^{sp})$ is specified, and its power injection is given in (21).

$$\begin{aligned}
(p_{i,Z}(\boldsymbol{y}) + jq_{i,Z}(\boldsymbol{y})) &= (e_i + jf_i)\left(\frac{e_i + jf_i}{\mathrm{Re}(z_i^{sp}) + j\,\mathrm{Im}(z_i^{sp})}\right)^* \\
&= \frac{e_i + jf_i}{(z_i^{sp})^2}\left[\left(e_i\,\mathrm{Re}(z_i^{sp}) + f_i\,\mathrm{Im}(z_i^{sp})\right) - j\left(f_i\,\mathrm{Re}(z_i^{sp}) - e_i\,\mathrm{Im}(z_i^{sp})\right)\right] \\
&= \frac{1}{(z_i^{sp})^2}\left[e_i\left(e_i\,\mathrm{Re}(z_i^{sp}) + f_i\,\mathrm{Im}(z_i^{sp})\right) + f_i\left(f_i\,\mathrm{Re}(z_i^{sp}) - e_i\,\mathrm{Im}(z_i^{sp})\right)\right] \\
&\quad + j\frac{1}{(z_i^{sp})^2}\left[f_i\left(e_i\,\mathrm{Re}(z_i^{sp}) + f_i\,\mathrm{Im}(z_i^{sp})\right) - e_i\left(f_i\,\mathrm{Re}(z_i^{sp}) - e_i\,\mathrm{Im}(z_i^{sp})\right)\right] \\
&= \underbrace{\frac{\mathrm{Re}(z_i^{sp})}{(z_i^{sp})^2}(e_i^2 + f_i^2)}_{p_{i,Z}(\boldsymbol{y})} + j\underbrace{\frac{\mathrm{Im}(z_i^{sp})}{(z_i^{sp})^2}(e_i^2 + f_i^2)}_{q_{i,Z}(\boldsymbol{y})}
\end{aligned} \tag{21}$$

For constant current component, the current $i_i^{sp} = \mathrm{Re}(i_i^{sp}) + j\,\mathrm{Im}(i_i^{sp})$ is specified, and its power injection is given in (22).

$$\begin{aligned}
(p_{i,I}(\boldsymbol{y}) + jq_{i,I}(\boldsymbol{y})) &= (e_i + jf_i)\left(\mathrm{Re}(i_i^{sp}) + j\,\mathrm{Im}(i_i^{sp})\right)^* \\
&= \left(e_i + jf_i\right)\left(\mathrm{Re}(i_i^{sp}) - j\,\mathrm{Im}(i_i^{sp})\right) \\
&= \underbrace{\left[e_i\,\mathrm{Re}(i_i^{sp}) + f_i\,\mathrm{Im}(i_i^{sp})\right]}_{p_{i,I}(\boldsymbol{y})} + j\underbrace{\left[f_i\,\mathrm{Re}(i_i^{sp}) - e_i\,\mathrm{Im}(i_i^{sp})\right]}_{q_{i,I}(\boldsymbol{y})}
\end{aligned} \tag{22}$$

For constant power component, the power $s_i^{sp} = \mathrm{Re}(s_i^{sp}) + j\,\mathrm{Im}(s_i^{sp}) = p_i^{sp} + jq_i^{sp}$ is specified, and its power injection is given in (23).

$$(p_{i,P}(\boldsymbol{y}) + jq_{i,P}(\boldsymbol{y})) = \underbrace{p_i^{sp}}_{p_{i,P}(\boldsymbol{y})} + j\underbrace{q_i^{sp}}_{q_{i,P}(\boldsymbol{y})} \tag{23}$$

Therefore, the total power injection is given by (24)-(25) where $\alpha, \beta$ are the percentages of each component (Z,I,P) in active power and reactive power injections, which satisfy the equalities $\alpha_{i,Z} + \alpha_{i,I} + \alpha_{i,P} = 1$ and $\beta_{i,Z} + \beta_{i,I} + \beta_{i,P} = 1$.

$$\begin{aligned}
p_{i,ZIP}(\boldsymbol{y}) &= \alpha_{i,Z}p_{i,Z}(\boldsymbol{y}) + \alpha_{i,I}p_{i,I}(\boldsymbol{y}) + \alpha_{i,P}p_{i,P}(\boldsymbol{y}) \\
&= \alpha_{i,Z}\frac{\mathrm{Re}(z_i^{sp})}{(z_i^{sp})^2}(e_i^2 + f_i^2) + \alpha_{i,I}\left[e_i\,\mathrm{Re}(i_i^{sp}) + f_i\,\mathrm{Im}(i_i^{sp})\right] + \alpha_{i,P}p_i^{sp}
\end{aligned} \tag{24}$$

$$\begin{aligned}
q_{i,ZIP}(\boldsymbol{y}) &= \beta_{i,Z}q_{i,Z}(\boldsymbol{y}) + \beta_{i,I}q_{i,I}(\boldsymbol{y}) + \beta_{i,P}q_{i,P}(\boldsymbol{y}) \\
&= \beta_{i,Z}\frac{\mathrm{Im}(z_i^{sp})}{(z_i^{sp})^2}(e_i^2 + f_i^2) + \beta_{i,I}\left[f_i\,\mathrm{Re}(i_i^{sp}) - e_i\,\mathrm{Im}(i_i^{sp})\right] + \beta_{i,P}q_i^{sp}
\end{aligned} \tag{25}$$

*2) DTs of Power Injections of ZIP Loads*

DTs of (24)-(25) are given by (26)-(27), respectively, using the transformation rules of DT method with details omitted.

$$\begin{aligned}
P_{i,ZIP}(\boldsymbol{Y}) &= \alpha_{i,Z}P_{i,Z}(\boldsymbol{Y}) + \alpha_{i,I}P_{i,I}(\boldsymbol{Y}) + \alpha_{i,P}P_{i,P}(\boldsymbol{Y}) \\
&= \alpha_{i,Z}\frac{\mathrm{Re}(z_i^{sp})}{(z_i^{sp})^2}(E_i(k)\otimes E_i(k) + F_i(k)\otimes F_i(k)) \\
&\quad + \alpha_{i,I}\left[E_i(k)\,\mathrm{Re}(i_i^{sp}) + F_i(k)\,\mathrm{Im}(i_i^{sp})\right] + \alpha_{i,P}p_i^{sp}\rho(k)
\end{aligned} \tag{26}$$

$$\begin{aligned}
Q_{i,ZIP}(\boldsymbol{Y}) &= \beta_{i,Z}Q_{i,Z}(\boldsymbol{Y}) + \beta_{i,I}Q_{i,I}(\boldsymbol{Y}) + \beta_{i,P}Q_{i,P}(\boldsymbol{Y}) \\
&= \beta_{i,Z}\frac{\mathrm{Im}(z_i^{sp})}{(z_i^{sp})^2}(E_i(k)\otimes E_i(k) + F_i(k)\otimes F_i(k)) \\
&\quad + \beta_{i,I}\left[F_i(k)\,\mathrm{Re}(i_i^{sp}) - E_i(k)\,\mathrm{Im}(i_i^{sp})\right] + \beta_{i,P}q_i^{sp}\rho(k)
\end{aligned} \tag{27}$$

**Proposition 3**: The transformed power injection equations of ZIP load model (26)-(27) respectively satisfy formally linear equations (28)-(29).

$$P_{i,ZIP}(\boldsymbol{Y}) = \boldsymbol{a}_{\mathrm{ZIP},i}\boldsymbol{Y}(k) + b_{ZIP,i} \tag{28}$$

$$Q_{i,ZIP}(\boldsymbol{Y}) = \boldsymbol{c}_{\mathrm{ZIP},i}\boldsymbol{Y}(k) + d_{ZIP,i} \tag{29}$$

Where $\boldsymbol{a}_{\mathrm{ZIP},i}$, $\boldsymbol{c}_{\mathrm{ZIP},i} \in \mathbb{R}^{1\times 2N}$ and $b_{\mathrm{ZIP},i}$, $d_{\mathrm{ZIP},i} \in \mathbb{R}$ are parameters given in the Appendix, together with the detailed proofs.

From Propositions 1-3, the following Propositions 4 and 5 are easily obtained, with detailed proofs omitted.

**Proposition 4**: Under the assumption of ZIP load, the nonlinear power flow equations (18)-(20) and (7)-(8) were converted to formally linear equations after DT. The detailed expressions of the linear equations are similar with those in Proposition 1, with slight modifications from (13)-(14) to (30)-(32).

$$0=\boldsymbol{a}_{\mathrm{P},i}\boldsymbol{Y}(k)-\Delta p_i\Lambda(k)+\varepsilon_i,\text{if } i\in\Omega_{PV} \tag{30}$$

$$\begin{gathered}0=\bar{\boldsymbol{a}}_{P,i}\boldsymbol{Y}(k)-\Delta p_i\Lambda(k)+\bar{\varepsilon}_i,\text{if } i\in\Omega_{PQ}\\ where\begin{cases}\bar{\boldsymbol{a}}_{P,i}=\boldsymbol{a}_{\mathrm{P},i}-\boldsymbol{a}_{\mathrm{ZIP},i}\\ \bar{\varepsilon}_i=\varepsilon_i+p_i^{sp}\rho(k)-b_{ZIP,i}\end{cases}\end{gathered} \tag{31}$$

$$\begin{gathered}0=\bar{\boldsymbol{a}}_{Q,i}\boldsymbol{Y}(k)-\Delta q_i\Lambda(k)+\bar{\mu}_i,\text{if } i\in\Omega_{PQ}\\ where\begin{cases}\bar{\boldsymbol{a}}_{Q,i}=\boldsymbol{a}_{\mathrm{Q},i}-\boldsymbol{c}_{\mathrm{ZIP},i}\\ \bar{\mu}_i=\mu_i+q_i^{sp}\rho(k)-d_{ZIP,i}\end{cases}\end{gathered} \tag{32}$$

**Proposition 5**: Under the assumption of ZIP load, the DTs (4b) of nonlinear power flow equations (4a) still satisfy the formally linear equation (4c). The detailed expressions of matrices $\boldsymbol{A}_{gy}$, $\mathbf{A}_{g\lambda}$, and $\boldsymbol{B}_g$ are similar with (17) in Proposition 2, with slight modifications on matrices $\boldsymbol{A}_{y,PQ}$ and $\boldsymbol{B}_{y,PQ}$ in (33).

$$\boldsymbol{A}_{y,\mathrm{PQ}}=\begin{bmatrix}\bar{\boldsymbol{a}}_{\mathrm{P},1}\\ \bar{\boldsymbol{a}}_{\mathrm{Q},1}\\ \vdots\\ \bar{\boldsymbol{a}}_{\mathrm{P},M}\\ \bar{\boldsymbol{a}}_{\mathrm{Q},M}\end{bmatrix}\quad \boldsymbol{B}_{\mathrm{PQ}}=\begin{bmatrix}\bar{\varepsilon}_1\\ \bar{\mu}_1\\ \vdots\\ \bar{\varepsilon}_M\\ \bar{\mu}_M\end{bmatrix} \tag{33}$$

## III. Conclusion

This letter extends the linear relationship in [2] to the nonlinear power flow equation with both constant power load and ZIP load model. It is proved that the nonlinear power flow equations satisfy formally linear equations after DT, under both load models. These results demonstrate the potential of the DT method to solve various power system nonlinear problems that are related with nonlinear power flow equations.

## Appendix

To make the proofs more compact, the following Lemma is first proved and then directly used as a conclusion to prove Proposition 1 and Proposition 3.

**Lemma:** For a quadratic nonlinear function $z(t)=x(t)y(t)$, its DT satisfies a formally linear equation in (34).

$$Z(k)=X(k)\otimes Y(k)=aX(k)+bY(k)+c \tag{34}$$

Especially, when $x(t)=y(t)$, (35) holds.

$$Z(k)=X(k)\otimes X(k)=2aX(k)+c \tag{35}$$

**Proof of Lemma:**

$$\begin{aligned}Z(k)=\ &X(k)\otimes Y(k)=\sum_{m=0}^{k}X(m)Y(k-m)\\ &=X(0)Y(k)+X(k)Y(0)+\sum_{m=1}^{k-1}X(m)Y(k-m)\end{aligned}$$

Therefore, (34) holds with *a, b* and *c* given below.

$$a=Y(0),b=X(0),c=\sum_{m=1}^{k-1}X(m)Y(k-m)$$

**Proof of Proposition 1:**

We only prove (13) here as an example, and the proofs of (14)-(16) are similar.

First, rewrite the RHS of (9) as

$$\begin{aligned}\mathrm{RHS}=&-\Delta p_i\Lambda(k)+\underbrace{g_{ii}\left(E_i(k)\otimes E_i(k)+F_i(k)\otimes F_i(k)\right)}_{\text{Term 1}}\\ &+\underbrace{\sum_{j=1,j\neq i}^{N}g_{ij}\left(E_i(k)\otimes E_j(k)+F_i(k)\otimes F_j(k)\right)}_{\text{Term 2}}\\ &+\underbrace{\sum_{j=1}^{N}b_{ij}\left(F_i(k)\otimes E_j(k)-E_i(k)\otimes F_j(k)\right)}_{\text{Term 3}}\end{aligned}$$

According to the Lemma, Term 1 is written as a formally linear equation about $E_i(k)$ and $F_i(k)$.

$$\begin{aligned}\text{Term 1}=&\ 2g_{ii}E_i(0)E_i(k)+2g_{ii}F_i(0)F_i(k)\\ &+g_{ii}\sum_{m=1}^{k-1}E_i(m)E_i(k-m)+g_{ii}\sum_{m=1}^{k-1}F_i(m)F_i(k-m)\end{aligned}$$

Then, both Term 2 and Term 3 are written as a formally linear equation about $E_i(k)$, $F_i(k)$, $E_j(k)$ and $F_j(k)$ below.

$$\begin{aligned}\text{Term 2}=&\sum_{\substack{j=1,\\ j\neq i}}^{N}g_{ij}\left(E_j(0)E_i(k)+E_i(0)E_j(k)\right)\\ &+\sum_{\substack{j=1\\ j\neq i}}^{N}g_{ij}\left(F_j(0)F_i(k)+F_i(0)F_j(k)\right)\\ &+\sum_{\substack{j=1\\ j\neq i}}^{N}g_{ij}\left(\sum_{m=1}^{k-1}E_i(m)E_j(k-m)+\sum_{m=1}^{k-1}F_i(m)F_j(k-m)\right)\end{aligned}$$

$$\begin{aligned}\text{Term 3}=&\sum_{j=1}^{N}b_{ij}\left(E_j(0)F_i(k)+F_i(0)E_j(k)\right)\\ &-\sum_{j=1}^{N}b_{ij}\left(F_j(0)E_i(k)+E_i(0)F_j(k)\right)\\ &+\sum_{j=1}^{N}b_{ij}\left(\sum_{m=1}^{k-1}F_i(m)E_j(k-m)-\sum_{m=1}^{k-1}E_i(m)F_j(k-m)\right)\end{aligned}$$

Finally, by summating the three terms and merging the similar terms, Equ. (13) holds with the vector $\boldsymbol{a}_{\mathrm{P},i}$ and parameter $\varepsilon_i$ given in (36) and (41). Similarly, Equ. (14)-(16) can be proved. Finally, vectors $\boldsymbol{a}_{\mathrm{P},i}$, $\boldsymbol{a}_{\mathrm{Q},i}$, $\boldsymbol{a}_{\mathrm{V},i}$, $\boldsymbol{a}_{\mathrm{E},i}$, $\boldsymbol{a}_{\mathrm{F},i}$, $\in\mathbb{R}^{1\times 2(N-1)}$ and parameters $\varepsilon_i$, $\mu_i$, $\zeta_i\in\mathbb{R}$ are in (36)-(43).

$$\begin{aligned} &\boldsymbol{a}_{\mathrm{P},i} = \begin{bmatrix}\alpha_{i1} & \beta_{i1} & \cdots & \alpha_{ij} & \beta_{ij} & \cdots\end{bmatrix}, \text{where} \\ &\alpha_{ij} = g_{ij}E_i(0) + b_{ij}F_i(0),\ \beta_{ij} = g_{ij}F_i(0) - b_{ij}E_i(0), \text{if } j \neq i \\ &\alpha_{ii} = \sum_{j=1}^{N}\left(g_{ij}E_j(0) - b_{ij}F_j(0)\right) + g_{ii}E_i(0) + b_{ii}F_i(0) \\ &\beta_{ii} = \sum_{j=1}^{N}\left(b_{ij}E_j(0) + g_{ij}F_j(0)\right) - b_{ii}E_i(0) + g_{ii}F_i(0) \end{aligned} \tag{36}$$

$$\begin{aligned} &\boldsymbol{a}_{\mathrm{Q},i} = \begin{bmatrix}\phi_{i1} & \psi_{i1} & \cdots & \phi_{ij} & \psi_{ij} & \cdots\end{bmatrix}, \text{where} \\ &\phi_{ij} = -b_{ij}E_i(0) + g_{ij}F_i(0),\ \psi_{ij} = -b_{ij}F_i(0) - g_{ij}E_i(0), \text{if } j \neq i \\ &\phi_{ii} = -\sum_{j=1}^{N}\left(b_{ij}E_j(0) + g_{ij}F_j(0)\right) - b_{ii}E_i(0) + g_{ii}F_i(0) \\ &\psi_{ii} = \sum_{j=1}^{N}\left(g_{ij}E_j(0) - b_{ij}F_j(0)\right) - g_{ii}E_i(0) - b_{ii}F_i(0) \end{aligned} \tag{37}$$

$$\boldsymbol{a}_{\mathrm{V},i} = \begin{bmatrix}0 & \cdots & 0 & 2E_i(0) & 2F_i(0) & 0 & \cdots & 0\end{bmatrix} \tag{38}$$

$$\boldsymbol{a}_{\mathrm{E},i} = \begin{bmatrix}0 & \cdots & 0 & E_i(0) & 0 & 0 & \cdots & 0\end{bmatrix} \tag{39}$$

$$\boldsymbol{a}_{\mathrm{F},i} = \begin{bmatrix}0 & \cdots & 0 & 0 & F_i(0) & 0 & \cdots & 0\end{bmatrix} \tag{40}$$

$$\begin{aligned} &\varepsilon_i = \sum_{j=1}^{N} g_{ij}c_{ij} + \sum_{j=1}^{N} b_{ij}d_{ij} - p_i\rho(k), \text{where} \\ &c_{ij} := \sum_{m=1}^{k-1} E_i(m)E_j(k-m) + \sum_{m=1}^{k-1} F_i(m)F_j(k-m) \\ &d_{ij} := \sum_{m=1}^{k-1} F_i(m)E_j(k-m) - \sum_{m=1}^{k-1} E_i(m)F_j(k-m) \end{aligned} \tag{41}$$

$$\mu_i = -\sum_{j=1}^{N} b_{ij}c_{ij} + \sum_{j=1}^{N} g_{ij}d_{ij} - q_i\rho(k) \tag{42}$$

$$\varsigma_i = c_{ii} - v_i^2\rho(k) \tag{43}$$

**Proof of Proposition 3:**

We only prove (28) here as an example, and the proof of (29) is similar. First, each component (Z, I, P) in (26) can be written as formally linear equations in (44)-(46).

$$\begin{aligned} P_{i,Z}(\boldsymbol{Y}) &= \frac{\mathrm{Re}(z_i^{sp})}{(z_i^{sp})^2}\left(E_i(k) + \sum_{m=1}^{k-1} E_i(m)E_i(k-m)\right) \\ &\quad + \frac{\mathrm{Re}(z_i^{sp})}{(z_i^{sp})^2}\left(2F_i(0)F_i(k) + \sum_{m=1}^{k-1} F_i(m)F_i(k-m)\right) \\ &= \left[\frac{2E_i(0)\,\mathrm{Re}(z_i^{sp})}{(z_i^{sp})^2}\right]E_i(k) + \left[\frac{2F_i(0)\,\mathrm{Re}(z_i^{sp})}{(z_i^{sp})^2}\right]F_i(k) \\ &\quad + \left[\frac{\mathrm{Re}(z_i^{sp})}{(z_i^{sp})^2}\left(\sum_{m=1}^{k-1} E_i(m)E_i(k-m) + \sum_{m=1}^{k-1} F_i(m)F_i(k-m)\right)\right] \\ &= \left[\frac{2E_i(0)\,\mathrm{Re}(z_i^{sp})}{(z_i^{sp})^2}\right]E_i(k) + \left[\frac{2F_i(0)\,\mathrm{Re}(z_i^{sp})}{(z_i^{sp})^2}\right]F_i(k) \\ &\quad + \left[\frac{\mathrm{Re}(z_i^{sp})}{(z_i^{sp})^2}c_{ii}\right] \end{aligned} \tag{44}$$

$$P_{i,I}(\boldsymbol{Y}) = \left[\mathrm{Re}(i_i^{sp})\right]E_i(k) + \left[\mathrm{Im}(i_i^{sp})\right]F_i(k) + \left[0\right] \tag{45}$$

$$P_{i,P}(\boldsymbol{Y}) = \left[0\right]E_i(k) + \left[0\right]F_i(k) + \left[p_i^{sp}\rho(k)\right] \tag{46}$$

Therefore, the total active power injection of ZIP load (26) is written as (47), which satisfy the linear equation (28) with vectors **a** and b given in (48)-(49).

$$\begin{aligned} P_{i,ZIP}(\boldsymbol{Y}) &= \alpha_{i,Z}P_{i,Z}(\boldsymbol{Y}) + \alpha_{i,I}P_{i,I}(\boldsymbol{Y}) + \alpha_{i,P}P_{i,P}(\boldsymbol{Y}) \\ &= \left[\alpha_{i,Z}\frac{2E_i(0)\,\mathrm{Re}(z_i^{sp})}{(z_i^{sp})^2} + \alpha_{i,I}\,\mathrm{Re}(i_i^{sp})\right]E_i(k) \\ &\quad + \left[\alpha_{i,Z}\frac{2F_i(0)\,\mathrm{Re}(z_i^{sp})}{(z_i^{sp})^2} + \alpha_{i,I}\,\mathrm{Im}(i_i^{sp})\right]F_i(k) \\ &\quad + \left[\alpha_{i,Z}\frac{\mathrm{Re}(z_i^{sp})}{(z_i^{sp})^2}c_{ii} + \alpha_{i,P}p_i^{sp}\rho(k)\right] \end{aligned} \tag{47}$$

$$\begin{aligned} &\boldsymbol{a}_{\mathrm{ZIP},i} = \begin{bmatrix}0 & \cdots & 0 & \tau_{1,i} & \tau_{2,i} & 0 & \cdots & 0\end{bmatrix} \\ &b_{ZIP,i} = \alpha_{i,Z}\frac{\mathrm{Re}(z_i^{sp})}{(z_i^{sp})^2}c_{ii} + \alpha_{i,P}p_i^{sp}\rho(k) \end{aligned} \tag{48}$$

$$\begin{aligned} &\tau_{1,i} = \alpha_{i,Z}\frac{2E_i(0)\,\mathrm{Re}(z_i^{sp})}{(z_i^{sp})^2} + \alpha_{i,I}\,\mathrm{Re}(i_i^{sp}) \\ &\tau_{2,i} = \alpha_{i,Z}\frac{2F_i(0)\,\mathrm{Re}(z_i^{sp})}{(z_i^{sp})^2} + \alpha_{i,I}\,\mathrm{Im}(i_i^{sp}) \end{aligned} \tag{49}$$

Similarly, the total reactive power injection of ZIP load (27) also satisfies the linear equation (29) with vectors **c** and d given in (50)-(51).

$$\begin{aligned} &\boldsymbol{c}_{\mathrm{ZIP},i} = \begin{bmatrix}0 & \cdots & 0 & \tau_{3,i} & \tau_{4,i} & 0 & \cdots & 0\end{bmatrix} \\ &d_{ZIP,i} = \beta_{i,Z}\frac{\mathrm{Im}(z_i^{sp})}{(z_i^{sp})^2}c_{ii} + \beta_{i,P}q_i^{sp}\rho(k) \end{aligned} \tag{50}$$

$$\begin{aligned} &\tau_{3,i} = \beta_{i,Z}\frac{2E_i(0)\,\mathrm{Im}(z_i^{sp})}{(z_i^{sp})^2} - \beta_{i,I}\,\mathrm{Im}(i_i^{sp}) \\ &\tau_{4,i} = \beta_{i,Z}\frac{2F_i(0)\,\mathrm{Im}(z_i^{sp})}{(z_i^{sp})^2} + \beta_{i,I}\,\mathrm{Re}(i_i^{sp}) \end{aligned} \tag{51}$$

## REFERENCES

[1] Y. Liu, K. Sun, R. Yao, B. Wang, "Power System Time Domain Simulation Using a Differential Transformation Method," in *IEEE Trans, Power Systems*, vol. 34, no. 5, pp. 3739-3748, Sept. 2019.

[2] Y. Liu, K. Sun, "Solving Power System Differential Algebraic Equations Using Differential Transformation," *IEEE Transactions on Power Systems*, in press (doi: 10.1109/TPWRS.2019.2945512)

[3] Y. Liu, K. Sun, " Differential Transformation of a Motor Load Model for Time-Domain Simulation," *arXiv preprint arXiv: 1908.09801*, 2019.

[4] P. Kundur, N. J. Balu, M. G. Lauby, *Power system stability and control*. McGraw-hill New York, 1994.

[5] F. Milano, *Power system modelling and scripting*. Springer Science & Business Media, 2010.

[6] R. D. Zimmerman, et al, “Matpower: Steady-state operations, planning and analysis tools for power systems research and education,” *IEEE Tran. Power Syst.,* vol. 26, no. 1, pp. 12–19, 2011.